\documentstyle[amssymb,12pt]{article}

\newtheorem{th}{Theorem}[section]
\newtheorem{lem}[th]{Lemma}
\newtheorem{prop}[th]{Proposition}

\newtheorem{defn}[th]{Definition}
\newenvironment{defn-new}{\begin{defn} \em}{\end{defn}}
\newtheorem{rem}[th]{Remark}
\newenvironment{rem-new}{\begin{rem} \em}{\end{rem}}
\newtheorem{ex}[th]{Example}
\newenvironment{ex-new}{\begin{ex} \em}{\end{ex}}
\newtheorem{prob}[th]{Problem}
\newenvironment{prob-new}{\begin{prob} \em}{\end{prob}}

\newenvironment{notation-new}{\begin{rem} \em}{\end{rem}}

\newenvironment{agr-new}{\begin{rem} \em}{\end{rem}}

\makeatletter \@addtoreset{equation}{section} \makeatother
\begin{document}

\begin{center}
{\Large {\bf On Lightlike Geometry of Indefinite Sasakian Statistical Manifolds}}%
\bigskip \bigskip

O\u{g}uzhan Bahad\i r
\end{center}

\bigskip \bigskip

\noindent {\bf Mathematics Subject Classification:} 53C15, 53C25, 53C40.

\noindent {\bf Keywords and phrases:} Lightlike hypersurface, Statistical
manifolds, Dual connections.

\medskip

\noindent {\bf Abstract.}
In this study, we introduce indefinite sasakian statistical manifolds and lightlike hypersurfaces of an indefinite sasakian statistical
manifold. Some relations among induced geometrical objects with respect to
dual connections in a lightlike hypersurface of an indefinite statistical manifold are
obtained. Some examples related to these concepts are also presented. Finally, we prove that an invariant lightlike submanifold of indefinite sasakian statistical manifold is an indefinite sasakian statistical manifold.

\section{Introduction\label{sect-intro}}

 Neural networks can be applied to solving numerous complex optimization problems in electromagnetic theory. Applied physicist B. Bartlett presented unsupervised machine learning model for computing approximate electromagnetic field solutions \cite{barlet}. In April 2019, the Event Horizon Telescope (EHT) collaboration released the first image of the shadow of a black hole with the help of deep learning algorithms. This image is direct evidence of the existence of black holes and general
theory of relativity \cite{Black}. This is also an indirect evidence of the existence of lightlike geometry in the universe.

 A statistical manifold, a new branch of mathematics,  is a generalization of the Riemannian manifold and is used to model the information;
  and also uses tools of differential geometry to study statistical inference, information loss and estimation
 \cite{Calin-Udr-2014-book}. Statistical manifolds also have many application areas such as neural networks, machine learning and artificial intelligence.

 Lightlike geometry is one of the most important research areas
in differential geometry and has many applications in physics and
mathematics, such as  general relativity, electromagnetism and black hole theory.

In 1975, Efron \cite{Efron-1975} first emphasized the
role of differential geometry in statistics. Differential geometrical tools were used by Amari to develop this idea
 \cite{Amari-1982}, \cite{Amari-1985-book}.
In 1989, Vos \cite{Vos-1989} obtained fundamental equations of geometry of submanifolds
of statistical manifolds. In 2009, hypersurfaces of a statistical manifold are studied by Furuhata \cite{Furuhata-2009}.
Many studies have been done on both statistical manifolds and lightlike geometry over the last few decades \cite{Aydin-MM-2015}-\cite{vilcu}. Hovewer,
no study combining these two notions has been done in the literature so far.

Motivated by these circumstances, in this study, we introduce the lightlike geometry of an indefinite sasakian statistical manifold. In Section~\ref{sect-prel}, we
present basic definitions and results about statistical manifolds and
lightlike hypersurfaces. In Section~\ref{sect-lightlike}, we show that
the induced connections on a lightlike hypersurface of a statistical manifold
need not be dual  and a lightlike hypersurface need not be a statistical
manifold. Moreover, we show that the second fundamental forms are not
degenerate. Finally, an example is given. In Section $4$, we introduce indefinite sasakian statistical manifolds and we obtain the characterization theorem of indefinite sasakian statistical manifolds. This section is concluded with two examples. In Section $5$, we consider lightlike hypersurfaces of indefinite sasakian statistical manifolds.
We characterize the parallelness, totaly geodeticity and integrability
of  some distributions. In this section we also give two examples. In Section $6$, we prove that an invariant lightlike submanifold of indefinite sasakian statistical manifold is an indefinite sasakian statistical manifold.

\section{Preliminaries\label{sect-prel}}

Let $(\overline{M},\bar{g})$ be an $(m+2)$-dimensional semi-Riemannian
manifold with ${\rm index}(\bar{g})=q\geq 1$. Let $(M,g)$ be a hypersurface
of $(\bar{M},\bar{g})$ with $g=\bar{g}|_{M}$. If the induced metric $g$ on $%
M $ is degenerate, then $M$ is called a lightlike (null or degenerate)
hypersurface (\cite{Duggal-Bejancu-1996}, \cite{Duggal-Jin-2007}, \cite%
{Duggal-Sahin-2010}). In this case, there exists a null vector field $\xi
\neq 0$ on $M$ such that
\begin{equation}
g\left( \xi ,X\right) =0,\qquad \forall \;X\in \Gamma \left( TM\right) .
\label{eq-null-1}
\end{equation}%
The radical or the null space of $T_{x}M$, at each point $x\in M$, is a
subspace $Rad~T_{x}M$ defined by
\begin{equation}
Rad~T_{x}M=\{\xi \in T_{x}M:g_{x}(\xi ,X)=0,\;X\in \Gamma (TM)\}.
\label{eq-null-2}
\end{equation}%
The dimension of $Rad~T_{x}M$ is called the nullity degree of $g$. We recall
that the nullity degree of $g$ for a lightlike hypersurface of $(\bar{M},%
\bar{g})$ is $1$. Since $g$ is degenerate and any null vector being
orthogonal to itself, $T_{x}M^{\perp }$ is also null and
\begin{equation}
Rad~T_{x}M=T_{x}M\cap T_{x}M^{\perp }.  \label{eq-null-3}
\end{equation}%
Since $\dim T_{x}M^{\perp }=1$ and $\dim Rad~T_{x}M=1,$ we have $%
Rad~T_{x}M=T_{x}M^{\perp }$. We call $Rad~TM$ a radical distribution and it
is spanned by the null vector field $\xi $. The complementary vector bundle $%
S(TM)$ of $Rad~TM$ in $TM$ is called the screen bundle of $M$. We note that
any screen bundle is non-degenerate. This means that
\begin{equation}
TM=Rad~TM\perp S(TM),  \label{eq-null-4}
\end{equation}%
with $\perp $ denoting the orthogonal-direct sum. The complementary vector
bundle $S(TM)^{\perp }$ of $S(TM)$ in $T\bar{M}$ is called screen
transversal bundle and it has rank $2$. Since $Rad~TM$ is a lightlike
subbundle of $S(TM)^{\perp }$ there exists a unique local section $N$ of $%
S(TM)^{\perp }$ such that
\begin{equation}
\bar{g}(N,N)=0,\quad \bar{g}(\xi ,N)=1.  \label{eq-null-5}
\end{equation}%
Note that $N$ is transversal to $M$ and $\{\xi ,N\}$ is a local frame field
of $S(TM)^{\perp }$ and there exists a line subbundle $ltr(TM)$ of $T\bar{M}$
, and it is called the lightlike transversal bundle, locally spanned by $N$.
Hence we have the following decomposition:
\begin{equation}
T\bar{M}=TM\oplus ltr(TM)=S(TM)\bot Rad~TM\oplus ltr(TM),  \label{eq-null-6}
\end{equation}%
where $\oplus $ is the direct sum but not orthogonal (\cite%
{Duggal-Bejancu-1996}, \cite{Duggal-Jin-2007}).
In view of the splitting $(\ref{eq-null-6})$, we have the following Gauss
and Weingarten formulas, respectively,
\begin{equation}
\bar{\nabla}_{X}Y=\nabla _{X}Y+h(X,Y),  \label{eq-null-7}
\end{equation}%
\begin{equation}
\bar{\nabla}_{X}N=-A_{N}X+\nabla _{X}^{t}N  \label{eq-null-8}
\end{equation}%
for any $X,Y\in \Gamma (TM)$, where $\nabla _{X}Y,\ A_{N}X\in \Gamma (TM)$
and $h(X,Y),\ \nabla _{X}^{t}N\in \Gamma (ltr(TM))$. If we set
\[
B(X,Y)=\bar{g}(h(X,Y),\xi )\quad {\rm and}\quad \tau (X)=\bar{g}(\nabla
_{X}^{t}N,\xi ),
\]%
then (\ref{eq-null-7}) and (\ref{eq-null-8}) become
\begin{equation}
\overline{\nabla }_{X}Y=\nabla _{X}Y+B(X,Y)N,  \label{eq-null-9}
\end{equation}%
\begin{equation}
\overline{\nabla }_{X}N=-A_{N}X+\tau (X)N,  \label{eq-null-10}
\end{equation}%
respectively. Here, $B$ and $A$ are called the second fundamental form and
the shape operator of the lightlike hypersurface $M$, respectively \cite%
{Duggal-Bejancu-1996}. Let $P$ be the projection of $S(TM)$ on $M$. Then,
for any $X\in \Gamma (TM)$, we can write
\begin{equation}
X=PX+\eta (X)\xi ,  \label{eq-null-11}
\end{equation}%
where $\eta $ is a $1$-form given by
\begin{equation}
\eta (X)=\bar{g}(X,N).  \label{eq-null-12}
\end{equation}%
From (\ref{eq-null-9}), we have
\begin{equation}
(\nabla _{X}g)(Y,Z)=B(X,Y)\eta (Z)+B(X,Z)\eta (Y),  \label{eq-null-13}
\end{equation}%
for all $X,Y,Z\in \Gamma (TM)$, where the induced connection $\nabla $ is a
non-metric connection on $M$. From (\ref{eq-null-4}), we have
\begin{equation}
\nabla _{X}W=\nabla _{X}^{\ast }W+h^{\ast }(X,W)=\nabla _{X}^{\ast
}W+C(X,W)\xi ,  \label{eq-null-14}
\end{equation}%
\begin{equation}
\nabla _{X}\xi =-A_{\xi }^{\ast }X-\tau (X)\xi  \label{eq-null-15}
\end{equation}%
for all $X\in \Gamma (TM)$, $W\in \Gamma (S(TM))$, where $\nabla _{X}^{\ast
}W$ and $A_{\xi }^{\ast }X$ belong to $\Gamma (S(TM))$. Here $C$, $A_{\xi
}^{\ast }$ and $\nabla ^{\ast }$ are called the local second fundamental
form, the local shape operator and the induced connection on $S(TM)$,
respectively. We also have
\begin{equation}
g(A_{\xi }^{\ast }X,W)=B(X,W),\ g(A_{\xi }^{\ast }X,N)=0,\ B(X,\xi )=0,\ \ \
g(A_{N}X,N)=0.  \label{eq-null-16}
\end{equation}%
Moreover, from the first and third equations of (\ref{eq-null-16}), we have
\begin{equation}
A_{\xi }^{\ast }\xi =0.  \label{eq-null-17}
\end{equation}%

Now we define some statistical basic concepts

\begin{defn-new}\label{sta}
{\rm \cite{Furuhata-2009}} Let $\widetilde{M}$ be a smooth manifold. Let $%
\widetilde{D}$ be an affine connection with the torsion tensor $T^{%
\widetilde{D}}$ and $\widetilde{g}$ a semi-Riemannian metric on $\widetilde{M%
}$. Then the pair $(\widetilde{D},\widetilde{g})$ is called a statistical
structure on $\widetilde{M}$ if
\begin{enumerate}
\item[{\rm (1)}] $(\widetilde{D}_{X}\widetilde{g})(Y,Z) - (\widetilde{D}_{Y}%
\widetilde{g})(X,Z) = \widetilde{g}(T^{\widetilde{D}}(X,Y),Z)$ \newline
for all $X,Y,Z\in \Gamma(T\widetilde{M})$, and

\item[{\rm (2)}] $T^{\widetilde{D}}=0$.
\end{enumerate}
\end{defn-new}

\begin{defn-new}
Let $(\widetilde{M},\widetilde{g})$ be a semi-Riemannian manifold. Two
affine connections $\widetilde{D}$ and $\widetilde{D}^{\ast }$ on $%
\widetilde{M}$ are said to be dual with respect to the metric $\widetilde{g}$%
, if
\begin{equation}
Z\widetilde{g}(X,Y)=\widetilde{g}(\widetilde{D}_{Z}X,Y)+\widetilde{g}(X,%
\widetilde{D}_{Z}^{\ast }Y)  \label{eq-dual-con}
\end{equation}%
for all $X,Y,Z\in \Gamma (T\widetilde{M})$.
\end{defn-new}

A statistical manifold will be represented by $(\widetilde{M},\widetilde{g},%
\widetilde{D},\widetilde{D}^{\ast })$. If $\widetilde{\nabla}$ is Levi-Civita
connection of $\widetilde{g}$, then
\begin{equation}
\widetilde{\nabla}=\frac{1}{2}(\widetilde{D}+\widetilde{D}^{\ast }).
\label{le}
\end{equation}%
In (\ref{eq-dual-con}), if we choose $\widetilde{D}^{\ast }=\widetilde{D}$
then Levi-Civita connection is obtained.
\begin{lem}
For statistical manifold $(\widetilde{M},\widetilde{g},%
\widetilde{D},\widetilde{D}^{\ast })$, we set
\begin{equation}
\overline{\mathbb{K}}=\widetilde{D}-\widetilde{\nabla}.\label{3}
\end{equation}
Then we have
\begin{equation}
\overline{\mathbb{K}}(X,Y)=\overline{\mathbb{K}}(Y,X),\;\widetilde{g}(\overline{\mathbb{K}}(X,Y),Z)=\widetilde{g}(\overline{\mathbb{K}}(X,Z),Y),\label{22}
\end{equation}
for any $X,Y,Z\in\Gamma(TM)$.

Conversely, for a Riemannian metric $g$, if $\overline{\mathbb{K}}$ satisfies (\ref{22}), the pair $(\widetilde{D}=\widetilde{\nabla}+\overline{\mathbb{K}},\widetilde{g})$ is a statistical structure on $\widetilde{M}$ {\rm\cite{Furuhata-Has-2017}}.
\end{lem}
Let $(M,g)$ be a submanifold of $(\widetilde{M},\widetilde{g})$. If $%
(M,g,D,D^{\ast })$ is a statistical manifold, then $(M,g,D,D^{\ast })$ is
called a statistical submanifold of $(\widetilde{M},\widetilde{g},\widetilde{%
D},\widetilde{D}^{\ast })$, where $D$, $D^{\ast }$ are affine dual
connections on $M$ and $\widetilde{D}$, $\widetilde{D}^{\ast }$ are affine
dual connections on $\widetilde{M}$ (see \cite{Amari-1985-book}, \cite%
{Furuhata-2009},\cite{Vos-1989}).

\section{Lightlike hypersurface of a statistical manifold\label%
{sect-lightlike}}

Let $(M,g)$ be a lightlike hypersurface of a statistical manifold $(%
\widetilde{M},\widetilde{g},\widetilde{D},\widetilde{D}^{\ast })$. Then,
Gauss and Weingarten formulas with respect to dual connections are given by
\cite{Furuhata-2009}
\begin{equation}
\widetilde{D}_{X}Y=D_{X}Y+B(X,Y)N,  \label{eq-D-tilde-X-Y}
\end{equation}%
\begin{equation}
\widetilde{D}_{X}N= -A_{N}X+\tau (X)N \label{eq-D-tilde-X-N}
\end{equation}%
\begin{equation}
\widetilde{D}_{X}^{\ast }Y=D_{X}^{\ast }Y+B^{\ast }(X,Y)N,
\label{eq-D-tilde-star-X-Y}
\end{equation}%
\begin{equation}
\widetilde{D}_{X}^{\ast }N=-A_{N}^{\ast }X+\tau ^{\ast }(X)N,  \label{eq-D-tilde-star-X-N}
\end{equation}%
for all $X,Y\in \Gamma (TM),\;N\in \Gamma (ltrTM)$, where $D_{X}Y$,$\
D_{X}^{\ast }Y$, $A_{N}X$,$\ A_{N}^{\ast }X\in \Gamma (TM)$ and
\[
B(X,Y)=\widetilde{g}(\widetilde{D}_{X}Y,\xi ),\quad {\tau }(X)=%
\widetilde{g}(\widetilde{D}_{X}N,\xi ),
\]%
\[
B^{\ast }(X,Y)=\widetilde{g}(\widetilde{D}_{X}^{\ast }Y,\xi ),\quad {\tau }^{\ast }%
(X)=\widetilde{g}(\widetilde{D}_{X}^{\ast }N,\xi ).
\]%
Here, $D$, $D^{\ast }$, $B$, $B^{\ast }$, ${A}_{N}$ and $A_{N}^{\ast }$ are
called the induced connections on $M$, the second fundamental forms and the
Weingarten mappings with respect to $\widetilde{D}$ and $\widetilde{D}^{\ast
}$, respectively. Using Gauss formulas and the equation (\ref{eq-dual-con}),
we obtain
\begin{eqnarray}
Xg(Y,Z) &=&g(\widetilde{D}_{X}Y,Z)+g(Y,\widetilde{D}_{X}^{\ast }Z),\
\nonumber \\
&=&g(D_{X}Y,Z)+g(Y,D_{X}^{\ast }Z)+B(X,Y)\eta (Z)+B^{\ast }(X,Z)\eta (Y).
\label{eq-LH-SM-1}
\end{eqnarray}

From the equation (\ref{eq-LH-SM-1}), we have the following result.

\begin{th} {\rm \cite{Oguz}}
\label{th-LH-SM-1} Let $(M,g)$ be a lightlike hypersurface of a statistical
manifold $(\widetilde{M},\widetilde{g},\widetilde{D},\widetilde{D}^{\ast })$%
. Then\/{\rm :}

\begin{enumerate}
\item[{\bf (i)}] Induced connections $D$ and $D^{\ast }$ need not be dual
.

\item[{\bf (ii)}] A lightlike hypersurface of a statistical manifold need
not be a statistical manifold with respect to the dual connections.
\end{enumerate}
\end{th}

Using Gauss and Weingarten formulas in (\ref{eq-LH-SM-1}), we get
\begin{eqnarray}
(D_{X}g)(Y,Z) + (D_{X}^{\ast }g)(Y,Z) &=& B(X,Y)\eta (Z) + B(X,Z)\eta(Y)
\nonumber \\
&&+B^{\ast }(X,Y)\eta (Z)+B^{\ast }(X,Z)\eta (Y).  \label{th-LH-SM-2}
\end{eqnarray}

\begin{prop} {\rm \cite{Oguz}}
\label{pro5} Let $(M,g)$ be a lightlike hypersurface of a statistical
manifold $(\widetilde{M},\widetilde{g},\widetilde{D},\widetilde{D}^{\ast })$%
. Then the following assertions are true\/{\rm :}

\begin{enumerate}
\item[{\bf (i)}] Induced connections $D$ and $D^{\ast }$ are symmetric
connection.

\item[{\bf (ii)}] The second fundamental forms $B$ and $B^{\ast }$ are
symmetric.
\end{enumerate}
\end{prop}

\noindent {\bf Proof.} We know that $T^{\widetilde{D}}=0$. Moreover,
\begin{eqnarray}
T^{\widetilde{D}}(X,Y)&=&\widetilde{D}_{X}Y-\widetilde{D}_{Y}X-[X,Y]
\nonumber \\
&=&D_{X}Y-D_{Y}X-[X,Y]+B(X,Y)N-B(Y,X)N=0.  \label{th-LH-SM-3}
\end{eqnarray}
Comparing the tangent and transversal components of (\ref{th-LH-SM-3}), we
obtain
\[
B(X,Y)=B(Y,X), \qquad T^{D}=0,
\]
where $T^{D}$ is the torsion tensor field of $D$. Thus, second fundamental
form $B$ is symmetric and induced connection $D$ is symmetric connection.

Similarly, it can be shown that the second fundamental form $B^{\ast}$ is
symmetric and the induced connection $D^{\ast}$ is a symmetric connection. $%
\blacksquare$

\medskip
Let $P$ denote the projection morphism of $\Gamma (TM)$ on $\Gamma (S(TM))$
with respect to the decomposition (\ref{eq-null-4}). Then, we have
\begin{equation}
D_{X}PY=\nabla _{X}PY+\overline{h}(X,PY),
\end{equation}%
\begin{equation}
D_{X}\xi =-\overline{A}_{\xi }X+\overline{\nabla }_{X}^{t}\xi =0
\end{equation}%
for all $X,Y\in \Gamma (TM)$ and $\xi \in \Gamma (RadTM)$, where $\nabla
_{X}PY$ and $\overline{A}_{\xi }X$ belong to $\Gamma (S(TM))$, $\nabla $ and
$\overline{\nabla }^{t}$ are linear connections on $\Gamma (S(TM))$ and $%
\Gamma (RadTM)$ respectively. Here, $\overline{h}$ and $\overline{A}$ are
called screen second fundamental form and screen shape operator of $S(TM)$,
respectively. If we define
\begin{equation}
C(X,PY)=g(\overline{h}(X,PY),N),
\end{equation}%
\begin{equation}
\varepsilon (X)=g(\overline{\nabla }_{X}^{t}\xi ,N),\;\forall X,Y\in \Gamma
(TM).
\end{equation}%
One can show that
\[
\varepsilon (X)=-\tau (X).
\]%
Therefore, we have
\begin{equation}
D_{X}PY=\nabla _{X}PY+C(X,PY)\xi ,  \label{c}
\end{equation}%
\begin{equation}
D_{X}\xi =-\overline{A}_{\xi }X-\tau (X)\xi =0,\;\forall X,Y\in \Gamma (TM).
\label{31}
\end{equation}%
Here $C(X,PY)$ is called the local screen fundamental form of $S(TM)$.

\medskip

Similarly, the relations of induced dual objects on $S(TM)$ are given by
\begin{equation}
D_{X}^{\ast }PY=\nabla _{X}^{\ast }PY+C^{\ast }(X,PY)\xi ,  \label{d}
\end{equation}%
\begin{equation}
D_{X}^{\ast }\xi =-\overline{A}_{\xi }^{\ast }X-\tau ^{\ast }(X)\xi
=0,\;\forall X,Y\in \Gamma (TM).  \label{33}
\end{equation}%
Using (\ref{eq-LH-SM-1}), (\ref{c}), (\ref{d}) and Gauss-Weingarten
formulas, the relationship between induced geometric objects are given by
\begin{equation}
B(X,\xi )+B^{\ast }(X,\xi )=0,\;g(A_{N}X+A_{N}^{\ast }X,N)=0,  \label{34}
\end{equation}%
\begin{equation}
C(X,PY)=g(A_{N}^{*}X,PY),\;C^{\ast }(X,PY)=g(A_{N}X,PY).  \label{18}
\end{equation}

Now, using the equation (\ref{34}) we can state the following result.

\begin{prop} {\rm \cite{Oguz}}
\label{pr6} Let $(M,g)$ be a lightlike hypersurface of a statistical
manifold $(\widetilde{M},\widetilde{g},\widetilde{D},\widetilde{D}^{\ast })$%
. Then second fundamental forms $B$ and $B^{\ast }$ are not degenerate.
\end{prop}

Additionally, due to $\widetilde{D}$ and $\widetilde{D}^{\ast }$ are dual
connections we obtain
\begin{equation}
B(X,Y)=g(\overline{A}_{\xi }^{\ast }X,Y)+B^{\ast }(X,\xi )\eta(Y),  \label{7}
\end{equation}%
\begin{equation}
B^{\ast }(X,Y)=g(\overline{A}_{\xi }X,Y)+B(X,\xi )\eta(Y).  \label{8}
\end{equation}%
Using (\ref{7}) and (\ref{8}) we get
\[
\overline{A}_{\xi }^{\ast }\xi +\overline{A}_{\xi }\xi =0.
\]

\begin{ex-new}
Let $(R_{2}^{4},\widetilde{g})$ be a $4$-dimensional semi-Euclidean space
with signature \allowbreak $(-,-,+,+)$ of the canonical basis $%
(\partial_{0},\ldots ,\partial_{3})$. Consider a hypersurface $M$ of $%
R_{2}^{4}$ given by
\[
x_{0}=x_{1}+\sqrt{2}\sqrt{x_{2}^{2}+x_{3}^{2}}.
\]
For simplicity, we set $f=\sqrt{x_{2}^{2}+x_{3}^{2}}$. It is easy to check
that $M$ is a lightlike hypersurface whose radical distribution $RadTM$ is
spanned by
\[
\xi=f(\partial_{0}-\partial_{1})+\sqrt{2}(x_{2}\partial_{2}+x_{3}%
\partial_{3}).
\]
Then the lightlike transversal vector bundle is given by
\[
ltr(TM)=Span\left\{{N=\frac{1}{4f^{2}}\left\{f(-\partial_{0}+\partial_{1})+\sqrt{2}%
(x_{2}\partial_{2}+x_{3}\partial_{3})\right\}}\right\}.
\]
It follows that the corresponding screen distribution $S(TM)$ is spanned by
\[
\{W_{1}=\partial_{0}+\partial_{1},\;W_{2}=-x_{3}\partial_{2}+x_{2}%
\partial_{3}\}.
\]
Then, by direct calculations we obtain
\[
\widetilde{\nabla}_{X}W_{1}=\widetilde{\nabla}_{W_{1}}X=0,
\]
\[
\widetilde{\nabla}_{W_{2}}W_{2}=-x_{2}\partial_{2}-x_{3}\partial_{3},
\]
\[
\widetilde{\nabla}_{\xi}\xi=\sqrt{2}\xi,\;\widetilde{\nabla}_{W_{2}}\xi=%
\widetilde{\nabla}_{\xi}W_{2}=\sqrt{2}W_{2},
\]
for any $X\in\Gamma(TM)$ \cite{Duggal-Sahin-2010}.

We define an affine connection $\widetilde{D}$ as follows
\begin{eqnarray}
&&\widetilde{D}_{X}W_{1}=\widetilde{D}_{W_{1}}X=0,\;\widetilde{D}%
_{W_{2}}W_{2}=-2x_{2}\partial_{2}  \nonumber \\
&&\widetilde{D}_{\xi}\xi=\sqrt{2}\xi,  \label{denk1} \\
&& \widetilde{D}_{W_{2}}\xi=\widetilde{D}_{\xi}W_{2}=\sqrt{2}W_{2}.  \nonumber
\end{eqnarray}
Then using (\ref{le}) we obtain
\begin{eqnarray}
&&\widetilde{D}^{\ast }_{X}W_{1}=\widetilde{D}^{\ast }_{W_{1}}X=0,\;%
\widetilde{D}^{\ast }_{W_{2}}W_{2}=-2x_{3}\partial_{3}  \nonumber \\
&&\widetilde{D}^{\ast }_{\xi}\xi=\sqrt{2}\xi,  \label{denk2} \\
&&\widetilde{D}^{\ast }_{W_{2}}\xi=\widetilde{D}^{\ast }_{\xi}W_{2}=\sqrt{2}%
W_{2}.  \nonumber
\end{eqnarray}
Then $\widetilde{D}$ and $\widetilde{D}^{*}$ are dual connections. Here, one can easily see that $T^{\widetilde{D}}=0$ and $\widetilde{D}\widetilde{g}=0$. From Definition \ref{sta}, we say that $(R_{2}^{4},\widetilde{g},\widetilde{D},\widetilde{D}^{\ast })$ is a
statistical manifold.

\end{ex-new}

\section{Indefinite sasakian statistical manifolds}

In order to call a differentiable semi-Riemannian manifold $(\widetilde{M},\widetilde{g})$ of dimension $n=2m+1$ as practically contact metric one, a $(1,1)$ tensor field $\widetilde{\varphi}$, a contravariant vector field $\nu$,  a $1-$ form $\eta$ and a Riemannian metric $\widetilde{g}$ should be admitted, which satisfy
\begin{eqnarray}
\widetilde{\varphi} \nu=0,\;\eta(\widetilde{\varphi} X)=0, \;\eta(\nu)=\epsilon,\label{2.1}\\
\widetilde{\varphi}^{2}(X)=-X+\eta(X)\nu,\; \widetilde{g}(X,\nu)=\epsilon\eta(X),\label{2.2}\\
\widetilde{g}(\widetilde{\varphi} X,\widetilde{\varphi} Y)=\widetilde{g}(X,Y)-\epsilon \eta(X)\eta(Y), \;\epsilon=\mp1\label{2.3}
\end{eqnarray}
for all the vector fields $X$, $Y$ on $\widetilde{M}$.  When a practically contact metric manifold performs
\begin{eqnarray}
(\widetilde{\nabla}_{X}\widetilde{\varphi})Y&=&\widetilde{g}(X,Y)\nu-\epsilon\eta(Y) X, \label{2.9}\\
\widetilde{\nabla}_{X}\nu&=&-\widetilde{\varphi}X,
\end{eqnarray}
$\widetilde{M}$ is regarded as an indefinite sasakian manifold. In this study, we assume that the vector field $\nu$ is spacelike.

\begin{defn}\label{ta}
 Let $ (\widetilde{g}, \widetilde{\varphi}, \nu)$ be an indefinite sasakian structure on $\widetilde{M}$. A quadruplet $(\widetilde{D} = \widetilde{\nabla} + \overline{\mathbb{K}}, \widetilde{g}, \widetilde{\varphi}, \nu)$ is called a indefinite sasakian statistical structure on $\widetilde{M}$ if $(\widetilde{D}, \widetilde{g})$ is a statistical structure on $\widetilde{M}$ and the formula
\begin{eqnarray}
\overline{\mathbb{K}}(X, \widetilde{\varphi} Y) = - \widetilde{\varphi} \overline{\mathbb{K}}(X, Y) \label{k}
\end{eqnarray}
holds for any $X, Y \in \Gamma(T\widetilde{M})$. Then $(\widetilde{M}, \widetilde{D}, \widetilde{g}, \widetilde{\varphi}, \nu)$ is said to an indefinite sasakian statistical manifold.
\end{defn}
An indefinite sasakian statistical manifold will be represented by $(\widetilde{M}, \widetilde{D}, \widetilde{g}, \widetilde{\varphi}, \nu)$. We remark that if $(\widetilde{M}, \widetilde{D}, \widetilde{g}, \widetilde{\varphi}, \nu)$ is an indefinite sasakian statistical manifold, so is $(\widetilde{M}, \widetilde{D}^{\ast }, \widetilde{g}, \varphi, \nu)$ \cite{Furuhata-Has-2016}, \cite{Furuhata-Has-2017}.

\begin{th} \label{th4.2}
Let $(\widetilde{M}, \widetilde{D}, \widetilde{g})$ be a statistical manifold and $(\widetilde{g}, \widetilde{\varphi}, \nu)$ an almost contact metric structure on $\widetilde{M}$. $(\widetilde{D}, \widetilde{g}, \widetilde{\varphi}, \nu)$ is an indefinite sasakian statistical struture if and only if the following conditions hold:
\begin{eqnarray}
&&\widetilde{D}_{X}{\varphi}Y-\widetilde{\varphi}\widetilde{D}^{\ast }_{X}Y=g(Y,X)\nu-\widetilde{g}(Y,\nu)X, \label{th1}\\
&&\widetilde{D}_{X}\nu=-\widetilde{\varphi}X+g(\widetilde{D}_{X}\nu,\nu)\nu \label{th2},
\end{eqnarray}
for all the vector fields $X$, $Y$ on $\widetilde{M}$.

\noindent{\bf Proof.}
Using (\ref{3}) we get
\begin{eqnarray}
\widetilde{D}_{X}\widetilde{\varphi}Y-\widetilde{\varphi}\widetilde{D}^{*}_{X}Y=(\widetilde{\nabla}_{X}\widetilde{\varphi})Y+\overline{\mathbb{K}}(X,\widetilde{\varphi}Y)+\widetilde{\varphi}\overline{\mathbb{K}}(X,Y)
\end{eqnarray}
for all the vector fields $X$, $Y$ on $\widetilde{M}$. If we consider Definition \ref{ta} and the equation (\ref{2.9}), we have the formula (\ref{th1}).
If we write $\widetilde{D}^{*}$ instead of $\widetilde{D}$ in (\ref{th1}), we have
\begin{eqnarray}
\widetilde{D}^{*}_{X}\widetilde{\varphi}Y-\widetilde{\varphi}\widetilde{D}_{X}Y=g(Y,X)\nu-\widetilde{g}(Y,\nu)X, \label{th3}
\end{eqnarray}
Substituting $\nu$ for $Y$ in (\ref{th3}), we have the equation (\ref{th2}).

Conversely using (\ref{th1}), we obtain
$$\widetilde{\varphi}\{\widetilde{D}_{X}{\varphi}^{2}Y-\widetilde{\varphi}\widetilde{D}^{\ast }_{X}\widetilde{\varphi}Y\}=0.$$
Assume (\ref{2.2}) and (\ref{th2}) as well, we get
$$0=-\widetilde{\varphi}\widetilde{D}_{X}Y+\widetilde{g}(Y,\nu)X-\widetilde{g}(X,\nu)\widetilde{g}(Y,\nu)\nu+\widetilde{D}^{*}_{X}\widetilde{\varphi}Y-\widetilde{g}(\widetilde{\varphi}X,\widetilde{\varphi}Y)\nu,$$
From (\ref{2.3}), this equation gives us (\ref{th3}).

Now, we will prove that (\ref{2.9}) and (\ref{k}) by using (\ref{th1}) and (\ref{th3}). Using (\ref{th1}) and (\ref{th3}), respectively, we have the following equations
$$(\widetilde{\nabla}_{X}\widetilde{\varphi})Y-g(Y,X)\nu+\widetilde{g}(Y,\nu)X=\overline{\mathbb{K}}(X,\widetilde{\varphi}Y)+\widetilde{\varphi}\overline{\mathbb{K}}(X,Y),$$
and
$$(\widetilde{\nabla}_{X}\widetilde{\varphi})Y-g(Y,X)\nu+\widetilde{g}(Y,\nu)X=-\overline{\mathbb{K}}(X,\widetilde{\varphi}Y)-\widetilde{\varphi}\overline{\mathbb{K}}(X,Y).$$
This last two equations verifies (\ref{2.9}) and (\ref{k}).

\end{th}

\begin{ex}\label{ex}
Let $\widetilde{M}=(R_{2}^{5},\widetilde{g})$ be a semi-Euclidean space,  where $\widetilde{g}$ is of the signature $(-,+,-,+,+)$ with respect to canonical basis $\{\frac{\partial}{\partial x_{1}},\frac{\partial}{\partial x_{2}}, \frac{\partial}{\partial y_{1}}, \frac{\partial}{\partial y_{2}}, \frac{\partial}{\partial z}\}$. Defining
\begin{eqnarray}
&&\eta=dz,\;\nu=\frac{\partial}{\partial z},\nonumber\\
&&\widetilde{\varphi}\left(\frac{\partial}{\partial x_{i}}\right)=-\frac{\partial}{\partial y_{i}},\;\widetilde{\varphi}\left(\frac{\partial}{\partial y_{i}}\right)=\frac{\partial}{\partial x_{i}},\;\widetilde{\varphi}\left(\frac{\partial}{\partial z}\right)=0, \nonumber
\end{eqnarray}
where $i=1,2$. It can easily see that $(\widetilde{\varphi},\nu,\eta,\widetilde{g})$  is an indefinite Sasakian structure on $R_{2}^{5}$.
If we choose $\overline{\mathbb{K}}(X,Y)=\widetilde{g}(Y,\nu)\widetilde{g}(X,\nu)\nu$, then $(\widetilde{D} = \widetilde{\nabla} + \overline{\mathbb{K}}, \widetilde{g}, \widetilde{\varphi}, \nu)$ is an indefinite Sasakian statistical structure on $\widetilde{M}$.
\end{ex}
\begin{ex} \label{ex1}
In a $5-$ dimensional real number space $\widetilde{M}=R^{5}$,  let $\{x_{i},y_{i},z\}_{1\leq i\leq 2}$ be cartesian coordinates on $\widetilde{M}$ and $\{\frac{\partial}{\partial x_{i}},\frac{\partial}{\partial y_{i}},\frac{\partial}{\partial z}\}_{1\leq i\leq 2}$ be the natural
field of frames. If we define $1-$ form $\eta$, a vector field $\nu$ and a tensor field $\widetilde{\varphi}$  as follows:
\begin{small}
\begin{eqnarray}
&&\eta=dz-y_{1}dx_{1}-x_{1}dy_{1},\;\nu=\frac{\partial}{\partial z},\nonumber\\
&&\widetilde{\varphi}(\frac{\partial}{\partial x_{1}})=-\frac{\partial}{\partial x_{2}},\;\widetilde{\varphi}(\frac{\partial}{\partial x_{2}})=\frac{\partial}{\partial x_{1}}+y_{1}\frac{\partial}{\partial z},\;\widetilde{\varphi}(\frac{\partial}{\partial y_{1}})=-\frac{\partial}{\partial y_{2}}\;\widetilde{\varphi}(\frac{\partial}{\partial y_{2}})=\frac{\partial}{\partial y_{1}}+x_{1}\frac{\partial}{\partial z},\;\widetilde{\varphi}(\frac{\partial}{\partial z})=0, \nonumber \\
\end{eqnarray}
\end{small}
It is easy to check (\ref{2.1}) and (\ref{2.2}) and thus $(\widetilde{\varphi},\nu,\eta)$ is an almost contact structure on $R^{5}$.
 Now, we define metric $\widetilde{g}$ on $R^{5}$ by \begin{small}
 \begin{eqnarray}
 \widetilde{g}&=&(y_{1}^{2}-1)dx_{1}^{2}-dx_{2}^{2}+(x_{1}^{2}+1)dy_{1}^{2}+dy_{2}^{2}+dz^{2}-y_{1}dx_{1}\otimes dz-y_{1}dz\otimes dx_{1}\nonumber \\ &+&x_{1}y_{1}dx_{1}\otimes dy_{1}+x_{1}y_{1}dy_{1}\otimes dx_{1}-x_{1}dy_{1}\otimes dz-x_{1}dz\otimes dy_{1}
 \end{eqnarray}
 \end{small}
 with respect to the natural field of frames. Then we can easily see that $(\widetilde{\varphi},\nu,\eta,\widetilde{g})$ is an indefinite Sasakian structure on $R^{5}$.
We set the difference tensor field $\overline{\mathbb{K}}$ as
\begin{eqnarray*}
\overline{\mathbb{K}}(X, Y) = \lambda\widetilde{g}(Y, \nu) \widetilde{g}(X, \nu) \nu,
\end{eqnarray*}
where $\lambda\in C^{\infty}(\widetilde{M})$. Then, $(\widetilde{D} = \widetilde{\nabla} + \overline{\mathbb{K}}, \widetilde{g}, \widetilde{\varphi}, \nu)$ is an indefinite Sasakian statistical structure on $\widetilde{M}$.
\end{ex}

\section{Lightlike hypersurfaces of indefinite sasakian statistical manifolds}

\begin{defn}
 Let $(M,g,D,D^{\ast })$ be a hypersurface of indefinite Sasakian statistical manifold
 $(\widetilde{M}, \widetilde{D}, \widetilde{g}, \widetilde{\varphi}, \nu)$. The quadruplet $(M,g,D,D^{\ast })$ is called lightlike hypersurface of indefinite Sasakian statistical manifold $(\widetilde{M}, \widetilde{D}, \widetilde{g}, \widetilde{\varphi}, \nu)$ if the induced metric $g$ is degenerate.
\end{defn}

Let $(\widetilde{M}, \widetilde{D}, \widetilde{g}, \widetilde{\varphi}, \nu)$ be a $(2m+1)-$ dimensional Sasakian statistical manifold and $(M,g)$ be a  lightlike hypersurface of $\widetilde{M}$, such that the
structure vector field $\nu$ is tangent to $M$. For any $\xi\in\Gamma(RadTM)$ and $N\in \Gamma (ltr TM)$, in view of (\ref{2.1})-(\ref{2.3}), we have
\begin{eqnarray}
\widetilde{g}(\xi,\nu)=0, \;\widetilde{g}(N,\nu)=0,\label{5.1}\\
\widetilde{\varphi}^{2}\xi=-\xi,\;\widetilde{\varphi}^{2}N=-N.
\end{eqnarray}
Also using (\ref{eq-D-tilde-X-Y}) and (\ref{th2}) we obtain
\begin{eqnarray}
B(\xi,\nu)=0,\;B(\nu,\nu)=0, \label{b1}\\
B^{*}(\xi,\nu)=0,\;B^{*}(\nu,\nu)=0. \label{b2}
\end{eqnarray}
\begin{prop}\label{pr}
Let $(\widetilde{M}, \widetilde{D}, \widetilde{g}, \widetilde{\varphi}, \nu)$ be a $(2m+1)-$ dimensional Sasakian statistical manifold and $(M,g,D,D^{\ast })$ be its lightlike hypersurface such that the structure vector field $\nu$ is tangent to $M$. Then we have
\begin{eqnarray}
g(\widetilde{\varphi}\xi,\xi)&=&0, \label{b1}\\
g(\widetilde{\varphi} \xi,N)&=&-g(\xi,\widetilde{\varphi}N)=-g(A_{N}^{*}\xi,\nu), \label{b2}\\
g(\widetilde{\varphi}\xi,\widetilde{\varphi}N)&=&1, \label{b11}
\end{eqnarray}
where $\xi$ is a local section of $RadTM$ and $N$ is a local section of $ltrTM$.

\noindent {\bf Proof.}
Using (\ref{th2}) and (\ref{eq-D-tilde-X-Y}), we have
\begin{eqnarray*}
g(\widetilde{\varphi}\xi,\xi)&=&g(-\widetilde{D}_{\xi}\nu+g(\widetilde{D}_{\xi}\nu,\nu)\nu,\xi)\\
&=&g(-D_{\xi}\nu-B(\xi,\nu)N,\xi),\\
&=&0
\end{eqnarray*}
and
\begin{eqnarray*}
g(\widetilde{\varphi}\xi,N)&=&g(-\widetilde{D}_{\xi}\nu+g(\widetilde{D}_{\xi}\nu,\nu)\nu,N)\\
&=&g(\nu,\widetilde{D}_{\xi}^{*}N),\\
&=&-g(A_{N}^{*}\xi,\nu).
\end{eqnarray*}
From (\ref{2.3}) and (\ref{5.1}), we have (\ref{b11}).

\end{prop}
Proposition \ref{pr}  makes it possible to make the following decompositions:
\begin{eqnarray}
S(TM)=\{\widetilde{\varphi}RadTM\oplus \widetilde{\varphi}ltr(TM)\}\bot L_{0}\bot \langle\nu\rangle\label{5.8},
\end{eqnarray}
where $L_{0}$ is non-degenerate and $\widetilde{\varphi}-$ invariant distribution of rank $2m-4$ on $M$. If we denote the following distributions on $M$
\begin{eqnarray}
L=RadTM\bot\widetilde{\varphi}RadTM\bot L_{0},\;L^{'}=\widetilde{\varphi}ltr(TM)\label{5.9},
\end{eqnarray}
then $L$ is invariant and $L^{'}$ is anti-invariant distributions under $\widetilde{\varphi}$. Also we have
\begin{eqnarray}
TM=L\oplus L^{'}\bot\langle\nu\rangle.
\end{eqnarray}

Now, we consider two null vector field $U$ and $W$ and their $1-$ forms $u$ and $w$ as follows:
\begin{eqnarray}
&&U=-\widetilde{\varphi}N,\;  u(X)=\widetilde{g}(X,W),\label{5.11}\\
&&W=-\widetilde{\varphi}\xi,\; w(X)=\widetilde{g}(X,U).\label{5.12}
\end{eqnarray}
Then, for any $X\in\Gamma(T\widetilde{M})$, we have
\begin{eqnarray}
X=SX+u(X)U,
\end{eqnarray}
where $S$ projection morphism of $T\widetilde{M}$ on the distribution $L$. Applying $\widetilde{\varphi}$ to last equation, we obtain
\begin{eqnarray}
\widetilde{\varphi}X&=&\widetilde{\varphi}SX+u(X)\widetilde{\varphi}U,\nonumber \\
\widetilde{\varphi}X&=&\varphi X+u(X)N, \label{58}
\end{eqnarray}
where $\varphi$ is a tensor field of type $(1,1)$ defined on $M$ by $\varphi X=\widetilde{\varphi}SX$.

Again, we apply $\widetilde{\varphi}$ to (\ref{58}) and using (\ref{2.1})-(\ref{2.3}) we have
\begin{eqnarray*}
\widetilde{\varphi}^{2}X&=&\widetilde{\varphi}\varphi X+u(X)\widetilde{\varphi}N,\\
-X+g(X,\nu)\nu&=&\varphi^{2}X-u(X)U.
\end{eqnarray*}
which means that
\begin{eqnarray}
\varphi^{2}X=-X+g(X,\nu)\nu+u(X)U.\label{4.23}
\end{eqnarray}
Now applying $\varphi$ to the equation (\ref{4.23}) and since $\varphi U=0$, we have $\varphi^{3}+\varphi=0$ which gives that $\varphi$ is an $f$--structure.

\begin{defn-new}
 Let $(M,g,D,D^{\ast })$ be a hypersurface of indefinite Sasakian statistical manifold
 $(\widetilde{M}, \widetilde{D}, \widetilde{g}, \widetilde{\varphi}, \nu)$. The quadruplet $(M,g,D,D^{\ast })$ is called screen semi-invariant  lightlike hypersurface of indefinite Sasakian statistical manifold $(\widetilde{M}, \widetilde{D}, \widetilde{g}, \widetilde{\varphi}, \nu)$ if
 \begin{eqnarray*}
 \widetilde{\varphi}(ltrTM)\subset S(TM),\\
 \widetilde{\varphi}(RadTM)\subset S(TM).
\end{eqnarray*}
\end{defn-new}
We remark that a hypersurface of indefinite Sasakian statistical manifold is screen semi-invariant lightlike hypersurface.

\begin{ex}
Let us recall the example \ref{ex}, Suppose that $M$ is a hypersurface of $R_{2}^{5}$ defined by
\begin{eqnarray*}
x_{1}=y_{2},
\end{eqnarray*}
Then $RadTM$ and $ltr(TM)$ are spanned by $\xi=\frac{\partial}{\partial x_{1}}+\frac{\partial}{\partial y_{2}}$ and $N=\frac{1}{2}\{-\frac{\partial}{\partial x_{1}}+\frac{\partial}{\partial y_{2}}\}$, respectively. Applying $\widetilde{\varphi}$ to this vector fields, we have
\begin{eqnarray*}
\widetilde{\varphi}\xi=\frac{\partial}{\partial x_{2}}-\frac{\partial}{\partial y_{1}},\;\widetilde{\varphi}N=\frac{1}{2}\left\{\frac{\partial}{\partial x_{2}}+\frac{\partial}{\partial y_{1}}\right\}.
\end{eqnarray*}
Thus M is a screen semi invariant lightlike hyperfurface of indefinite sasakian statistical manifold $R_{2}^{5}$.
\end{ex}

\begin{ex}
Let M be a hypersurface of $(\widetilde{\varphi},\nu,\eta,\widetilde{g})$ on $\widetilde{M}=R^{5}$ in Example\ref{ex1}, Suppose that $M$ is a hypersurface of $R_{2}^{5}$ defined by
\begin{eqnarray*}
x_{2}=y_{2},
\end{eqnarray*}
Then the tangent space $TM$ is spanned by $\{U_{i}\}_{1\leq i\leq 4}$, where $U_{1}=\frac{\partial}{\partial x_{1}}$, $U_{2}=\frac{\partial}{\partial x_{2}}+\frac{\partial}{\partial y_{2}}$, $U_{3}=\frac{\partial}{\partial y_{1}}$, $U_{4}=\nu$.
$RadTM$ and $ltr(TM)$ are spanned by $\xi=U_{2}$ and $N=\frac{1}{2}\{-\frac{\partial}{\partial x_{2}}+\frac{\partial}{\partial y_{2}}\}$, respectively. Applying $\widetilde{\varphi}$ to this vector fields, we have
\begin{eqnarray*}
\widetilde{\varphi}\xi=U_{1}+U_{3}+(x_{1}+y_{1})U_{4},\;\widetilde{\varphi}N=\frac{1}{2}\left\{-U_{1}+U_{3}+(x_{1}-y_{1})U_{4}\right\}.
\end{eqnarray*}
Thus M is a screen semi invariant lightlike hyperfurface of indefinite sasakian statistical manifold $\widetilde{M}$.
\end{ex}
In view of (\ref{5.11}) and (\ref{5.12}), we have
\begin{eqnarray*}
\widetilde{g}(U,W)=1.
\end{eqnarray*}
Thus $\langle U\rangle \oplus \langle W\rangle$ is non-degenerate vector budle of $S(TM)$ with rank $2$. If we consider (\ref{5.8}) and (\ref{5.9}), we get
\begin{eqnarray}
S(TM)=\{U\oplus W\}\bot L_{0}\bot \langle\nu\rangle\label{5.33},
\end{eqnarray}
and
\begin{eqnarray}
L=RadTM\bot \langle W\rangle \bot L_{0},\;L^{'}=\langle U\rangle\label{5.34}.
\end{eqnarray}
Thus, for any $X\in\Gamma(TM)$, we can write
\begin{eqnarray}
X=PX+QX+g(X,\nu)\nu,\label{5.35}
\end{eqnarray}
where P and Q are projections of $TM$ into $L$ and $L^{'}$. Thus, we can write $QX=u(X)U$. Using  (\ref{2.1}),(\ref{2.2}), (\ref{2.3}), (\ref{58}) and (\ref{5.35}), we have
$$\varphi^{2}X=-X+g(X,\nu)\nu+u(X)U.$$
where $\widetilde{\varphi}PX=\varphi X$. We can easily see that
\begin{eqnarray}
g(\varphi X,\varphi Y)=g(X,Y)-g(X,\nu)g(Y,\nu)-u(X)w(Y)-u(Y)w(X),
\end{eqnarray}
for any $X,Y\in\Gamma(TM)$. Also we have the following identities:
\begin{eqnarray}
g(\varphi X, Y)=g(X,\varphi Y)-u(X)\eta(Y)-u(Y)\eta(X),\\
\varphi \nu=0,\;g(\varphi X,\nu)=0.
\end{eqnarray}
Thus, we have the following proposition
\begin{prop}
Let $(M,g,D,D^{\ast })$ be a lightlike hypersurface of indefinite Sasakian statistical manifold
 $(\widetilde{M}, \widetilde{D}, \widetilde{g}, \widetilde{\varphi}, \nu)$. Then $\varphi$ need not be a almost contact structure.
\end{prop}

\begin{lem} \label{lem52}
 Let $(M,g,D,D^{\ast })$ be a lightlike hypersurface of indefinite Sasakian statistical manifold  $(\widetilde{M}, \widetilde{D}, \widetilde{g}, \widetilde{\varphi}, \nu)$. For any $X,Y\in\Gamma(TM)$, we have the following identities:
\begin{eqnarray}
D_{X}\varphi Y-\varphi D^{*}_{X}Y=u(Y)A_{N}X-B^{*}(X,Y)U+g(X,Y)\nu-g(\nu,Y)X, \label{1}\\
D_{X}(u(Y))-u(D^{*}_{X}Y)=-B(X,\varphi Y)-u(Y)\tau(X) \label{2}
\end{eqnarray}
 \noindent{\bf Proof.}
Using Gauss and Weingarten formulas in (\ref{th1}) we obtain \begin{small}
\begin{eqnarray}
D_{X}\varphi Y+B(X,\widetilde{\varphi}Y)+D_{X}(u(Y))N-u(Y)A_{N}X+u(Y)\tau(X)N-\varphi\nabla^{*}_{X}Y+B^{*}(X,Y)U \nonumber\\
=g(X,Y)\nu-g(\nu,Y)X
\end{eqnarray}
If we take tangential and transversal parts of this last
equation we have (\ref{1}) and (\ref{2}).
\end{small}
\end{lem}
Similarly, we have the following lemma
\begin{lem}\label{lem53}
 Let $(M,g,D,D^{\ast })$ be a lightlike hypersurface of indefinite Sasakian statistical manifold  $(\widetilde{M}, \widetilde{D}, \widetilde{g}, \widetilde{\varphi}, \nu)$. For any $X,Y\in\Gamma(TM)$, we have the following identities:
\begin{eqnarray}
D^{*}_{X}\varphi Y-\varphi D_{X}Y=u(Y)A^{*}_{N}X-B(X,Y)U+g(X,Y)\nu-g(\nu,Y)X, \label{1*}\\
D^{*}_{X}(u(Y))-u(D_{X}Y)=-B^{*}(X,\varphi Y)-u(Y)\tau^{*}(X) \label{2*}
\end{eqnarray}
 \end{lem}
Lemme(\ref{lem52}) and Lemma(\ref{lem53}) are give us the following theorem.
\begin{th}
A lightlike hypersurface $M$ of an indefinite Sasakian statistical manifold $\widetilde{M}$ need
not be a statistical manifold.
\end{th}

\begin{prop}
 Let $(M,g,D,D^{\ast })$ be a lightlike hypersurface of indefinite Sasakian statistical manifold  $(\widetilde{M}, \widetilde{D}, \widetilde{g}, \widetilde{\varphi}, \nu)$. For any $X,Y\in\Gamma(TM)$, we have the following expressions:
  \\
 \textbf{(i)} If the vector field $U$ is parallel with respect to $\nabla^{*}$, then
 \begin{eqnarray}
A_{N}X=u(A_{N}X)U+\tau(A_{N}X)\nu.\;\tau(X)=0
\end{eqnarray}
 \textbf{(ii)} If the vector field $U$ is parallel with respect to $\nabla$, then
 \begin{eqnarray}
A_{N}^{*}X=u(A_{N}^{*}X)U+\tau(A_{N}^{*}X)\nu.\;\tau^{*}(X)=0\label{ii}
\end{eqnarray}
 \noindent{\bf Proof.}
 Replacing $Y$ in (\ref{1}) by $U$, we obtain $$-\varphi D^{*}_{X}Y=A_{N}X-B^{*}(X,U)U+g(X,U)\nu.$$
 Applying $\varphi$ to this equation and using (\ref{4.23}), we get
 $$D_{X}^{*}U-g(D_{X}^{*}U,\nu)\nu-u(D_{X}^{*}U)U=\varphi A_{N}X.$$
If $U$ is parallel with respect to $\nabla^{*}$ then $\varphi A_{N}X=0$. From (\ref{58}) we have $\widetilde{\varphi}(A_{N}X)=u(A_{N}X)N$. Applying $\widetilde{\varphi}$ this and using (\ref{2.2}) we obtain $A_{N}X=u(A_{N}X)U+\tau(A_{N}X)\nu$.
Also, if we write $U$ instead of $Y$ in the equation (\ref{2}), we have $\tau(X)=0$.

(\ref{ii}) can also be easily obtained by similar method.

\end{prop}

\begin{prop}
 Let $(M,g,D,D^{\ast })$ be a lightlike hypersurface of indefinite Sasakian statistical manifold  $(\widetilde{M}, \widetilde{D}, \widetilde{g}, \widetilde{\varphi}, \nu)$. For any $X,Y\in\Gamma(TM)$, we have the following expressions:
  \\
 \textbf{(i)} If the vector field $W$ is parallel with respect to $\nabla^{*}$, then
 \begin{eqnarray}
\overline{A}^{*}_{\xi}X=g( \overline{A}^{*}_{\xi}X,\nu)\nu+u( \overline{A}^{*}_{\xi}X)U,\; \tau^{*}(X)=0.\label{5.23}
\end{eqnarray}
 \textbf{(ii)} If the vector field $W$ is parallel with respect to $\nabla$, then
 \begin{eqnarray}
\overline{A}_{\xi}X=g( \overline{A}_{\xi}X,\nu)\nu+u( \overline{A}_{\xi}X)U,\; \tau(X)=0.\label{5.24}
\end{eqnarray}
 \noindent{\bf Proof.}
If we write $\xi$ instead of $Y$ in the equation (\ref{1}), we obtain
$$D_{X}\varphi \xi-\varphi D_{X}^{*}\xi=-B^{*}(X,\xi)U.$$
If $W$ is parallel with respect to $D$, using (\ref{33}) and (\ref{5.12}) in this equation, we obtain
$$\varphi \overline{A}^{*}_{\xi}X-\tau^{*}(X)W=-B^{*}(X,\xi)U.$$
Applying $\widetilde{\varphi}$ this and using (\ref{4.23}) we have
$$- \overline{A}^{*}_{\xi}X+g( \overline{A}^{*}_{\xi}X,\nu)\nu+u( \overline{A}^{*}_{\xi}X)U=\tau^{*}(X)\xi$$
If we take screen and radical parts of this last
equation we have (\ref{5.23}).

Similarly, we can easily see the equation (\ref{5.24}).
\end{prop}
\begin{defn-new}
{\rm (\cite{Furuhata-Has-2016}, \cite{Kurose-2002})} Let $(M,g)$ be a
hypersurface of a statistical manifold $(\widetilde{M},\widetilde{g},%
\widetilde{D},\widetilde{D}^{\ast })$.

\begin{enumerate}
\item[{\bf (i)}] $M$ is called \emph{totally geodesic with respect to} $\widetilde{D%
}$ if $B=0$.

\item[{\bf (ii)}] $M$ is called \emph{totally geodesic with respect to} $\widetilde{%
D}^{\ast }$ if $B^{\ast }=0$.

\end{enumerate}
\end{defn-new}
\begin{th}
Let $(M,g,D,D^{\ast })$ be a lightlike hypersurface of indefinite Sasakian statistical manifold  $(\widetilde{M}, \widetilde{D}, \widetilde{g}, \widetilde{\varphi}, \nu)$. \\
\textbf{(i)} $M$ is totally geodesic with respect to $\widetilde{D%
}$ if and only if
\begin{eqnarray}
D_{X}\varphi Y-\varphi D_{X}^{*}Y=g(X,Y)\nu,\;\forall X\in\Gamma(TM),\;Y\in\Gamma(L),\\
A_{N}X=-\varphi D_{X}^{*}U-g(X,U)\nu,\;\forall X\in\Gamma(TM).
\end{eqnarray}
\textbf{(ii)} $M$ is totally geodesic with respect to $\widetilde{D%
}^{*}$ if and only if
\begin{eqnarray}
D^{*}_{X}\varphi Y-\varphi D_{X}Y=g(X,Y)\nu,\;\forall X\in\Gamma(TM),\;Y\in\Gamma(L),\\
A^{*}_{N}X=-\varphi D_{X}U-g(X,U)\nu,\;\forall X\in\Gamma(TM).
\end{eqnarray}
\noindent{\bf Proof.}
For any $X\in\Gamma(TM)$ we know that $u(Y)=0$. Then the equations (\ref{1}) and (\ref{1*}) are reduced to the equations, respectively
 \begin{eqnarray}
 D_{X}\varphi Y-\varphi D^{*}_{X}Y=-B^{*}(X,Y)U+g(X,Y)\nu, \label{5.29}\\
 D^{*}_{X}\varphi Y-\varphi D_{X}Y=-B(X,Y)U+g(X,Y)\nu.\label{5.30}
\end{eqnarray}
On the other hand, replacing $Y$ by $U$ in (\ref{1}) and (\ref{1*}), respectively, we also have
 \begin{eqnarray}
 A_{N}X=-\varphi D_{X}^{*}U+B^{*}(X,U)U-g(X,U)\nu \label{5.31}\\
A^{*}_{N}X=-\varphi D_{X}U+B(X,U)U-g(X,U)\nu. \label{5.32}
\end{eqnarray}
If taking into account (\ref{5.29}), (\ref{5.30}), (\ref{5.31}) and (\ref{5.32}), we can easily obtain our assertion.
\end{th}

The following two theorems give a characterization of the integrability of distributions $L\bot\langle\nu\rangle$ and $L^{'}\bot\langle\nu\rangle$, respectively.
\begin{th}
Let $(M,g,D,D^{\ast })$ be a screen semi-invariant hypersurface of indefinite Sasakian statistical manifold
 $(\widetilde{M}, \widetilde{D}, \widetilde{g}, \widetilde{\varphi}, \nu)$. The following assertions are equivalent:\\
 \textbf{(i)} The distribution $L\bot\langle\nu\rangle$ is integrable. \\
 \textbf{(ii)} $B^{*}(X,\varphi Y)=B^{*}(\varphi X, Y)$, for all $X,Y\in\Gamma(L\bot\langle\nu\rangle)$,\\
 \textbf{(iii)}  $B(X,\varphi Y)=B(\varphi X, Y)$, for all $X,Y\in\Gamma(L\bot\langle\nu\rangle)$.

 \noindent{\bf Proof.}
 We know that $X\in\Gamma(L\bot\langle\nu\rangle)$ if and only if $u(X)=\widetilde{g}(X,W)=0$.
 For any $X,Y\in\Gamma(L\bot\langle\nu\rangle)$, using (\ref{eq-D-tilde-X-Y}) and (\ref{58}), we obtain
$$
u[X,Y]=-u(D_{X}Y)+u(D_{Y}X).
$$
From (\ref{2}), we have
$$
u[X,Y]=B^{*}(Y,\varphi X)-B^{*}(\varphi Y,X).
$$
This gives the equivalence between (i) and (ii).
 Similarly we can easily see that the relation (i) and (iii).
\end{th}

\begin{th}
Let $(M,g,D,D^{\ast })$ be a screen semi-invariant hypersurface of indefinite Sasakian statistical manifold
 $(\widetilde{M}, \widetilde{D}, \widetilde{g}, \widetilde{\varphi}, \nu)$. The following assertions are equivalent:\\
 \textbf{(i)} The distribution $L^{'}\bot\langle\nu\rangle$ is integrable.\\
\textbf{(ii)}  $A^{*}_{\widetilde{\varphi} X}Y-A^{*}_{\widetilde{\varphi} Y}X=g(X,\nu)Y-g(Y,\nu)X$, for all $X,Y\in\Gamma(L^{'}\bot\langle\nu\rangle)$.\\
 \textbf{(ii)}$A_{\widetilde{\varphi} X}Y-A_{\widetilde{\varphi} Y}X=g(X,\nu)Y-g(Y,\nu)X$, for all $X,Y\in\Gamma(L^{'}\bot\langle\nu\rangle)$.

  \noindent{\bf Proof.}
$X\in\Gamma(L^{'}\bot\langle\nu\rangle)$ if and only if $\varphi X=0$.
 For any $X,Y\in\Gamma(L\bot\langle\nu\rangle)$, using (\ref{eq-D-tilde-X-N}), (\ref{eq-D-tilde-star-X-Y}) and (\ref{58}) in (\ref{th1}), we have
$$
\varphi D^{*}_{X}Y=-g(X,Y)\nu+\widetilde{g}(Y,\nu)X-A_{\widetilde{\varphi}Y}X+B^{*}(X,Y)U.
$$
Therefore, we can get
$$
\varphi[X,Y]=-A_{\widetilde{\varphi}Y}X+A_{\widetilde{\varphi}X}Y+\widetilde{g}(Y,\nu)X-\widetilde{g}(X,\nu)Y.
$$
This gives the equivalence between (i) and (ii).
 Similarly we can easily see that the relation (i) and (iii).
\end{th}
\section{Invariant submanifolds}
Let $(M,g,D,D^{\ast })$ be an invariant lightlike submanifold of indefinite Sasakian statistical manifold
 $(\widetilde{M}, \widetilde{D}, \widetilde{g}, \widetilde{\varphi}, \nu)$. if $M$ is tangent to the structure vector field $\nu$, then $\nu$ belongs to $S(TM)$ (see \cite{Duggal-Sahin-2010}). For invariant lightlike submanifold, we have the following expressions:
 \begin{eqnarray}
 \widetilde{\varphi}(S(TM))=S(TM),\;\widetilde{\varphi}(RadTM)=RadTM
 \end{eqnarray}
\begin{prop}
Let $(M,g,D,D^{\ast })$ be an invariant lightlike submanifold of indefinite Sasakian statistical manifold
 $(\widetilde{M}, \widetilde{D}, \widetilde{g}, \widetilde{\varphi}, \nu)$ such that
the structure vector field $\nu$ is tangent to $M$. For any $X,Y\in\Gamma(TM)$, we have the following identities:
\begin{eqnarray}
D_{X}\varphi Y-\varphi D^{*}_{X}Y&=&g(X,Y)\nu-g(\nu,Y)X,\label{6.2}\\
h(X,\widetilde{\varphi} Y)&=&\widetilde{\varphi}h^{*}(X,Y)\label{6.3},
\end{eqnarray}
where $h$ and $h^{*}$ are second fundemental forms for affine dual connections $\widetilde{D}$ and $\widetilde{D}^{*}$, respectively.

\noindent{\bf Proof.}
Using (\ref{58}) and Gauss formula in (\ref{th1}), we obtain
$$D_{X}\varphi Y+h(X,\widetilde{\varphi} Y)-\varphi D^{*}_{X}Y-\widetilde{\varphi}h^{*}(X,Y)=g(X,Y)\nu-g(\nu,Y)X.$$
If we take tangential and transversal parts of this last equation, our claim is proven.
\end{prop}
Similar to the above proposition, the following proposition is given for dual connection $D^{*}$.
\begin{prop}
Let $(M,g,D,D^{\ast })$ be an invariant lightlike submanifold of indefinite Sasakian statistical manifold
 $(\widetilde{M}, \widetilde{D}, \widetilde{g}, \widetilde{\varphi}, \nu)$ such that
the structure vector field $\nu$ is tangent to $M$. For any $X,Y\in\Gamma(TM)$, we have the following identities:
\begin{eqnarray}
D^{*}_{X}\varphi Y-\varphi D_{X}Y&=&g(X,Y)\nu-g(\nu,Y)X,\\
h^{*}(X,\widetilde{\varphi} Y)&=&\widetilde{\varphi}h(X,Y)\label{6.5},
\end{eqnarray}
where $h$ and $h^{*}$ are second fundemental forms for affine dual connections $\widetilde{D}$ and $\widetilde{D}^{*}$, respectively.
\end{prop}
From the equations (\ref{6.3}) and (\ref{6.5}), we have
\begin{eqnarray}
h(X,\nu)=0,\;h^{*}(X,\nu)=0, \label{h}
\end{eqnarray}

A lightlike submanifold may not be an indefinite sasakian statistical manifold. The following theorem gives a case where this can happen.
\begin{th}
An invariant lightlike submanifold of indefinite Sasakian statistical manifold is an indefinite sasakian statistical manifold.

\noindent{\bf Proof.}
In a invariant lightlike submanifold, $u(X)=0$, for any $X\in\Gamma(TM)$. Then from (\ref{58}) we have
$$\varphi^{2}X=-X+g(X,\nu)\nu.$$
Since $\widetilde{\varphi}X=\varphi X$, using (\ref{2.1}), (\ref{2.2}) and (\ref{2.3}), we obtain
\begin{eqnarray}
\varphi \nu=0,\;\eta(\varphi X)=0, \\
\widetilde{g}(\varphi X,\varphi Y)=g(X,Y)-\eta(X)\eta(Y).
\end{eqnarray}
Then $(g, \varphi, \nu)$ is an almost contact metric structure.

Using (\ref{eq-LH-SM-1}), we get
\begin{eqnarray}
Xg(\varphi Y,\varphi Z)=g(D_{X}\varphi Y,\varphi Z)+g(\varphi Y,D_{X}^{\ast }\varphi Z).
\end{eqnarray}
This equation says that $D$ and $D^{*}$ are dual connections. Moreover torsion tensor of the connection $D$ is equal zero. Then, the equations (\ref{eq-LH-SM-1}) and (\ref{th-LH-SM-2}) tell us that $(D,g)$ is a statistical structure.

If we consider Gauss formula and (\ref{th2}) we obtain
\begin{eqnarray}
D_{X}\nu=-\varphi X+g(D_{X}\nu,\nu)\nu\label{6.10}.
\end{eqnarray}
If we consider (\ref{6.2}) and (\ref{6.10}) in the theorem \ref{th4.2}, our assertion are proven.
\end{th}
\begin{center}
\bigskip {\bf 4. Conclusion and future work}
\end{center}
In the present paper, firstly we have studied lightlike geometry of statistical manifolds, Later, we have introduced lightlike geometry of an indefinite sasakian statistical manifold which is a new classification of statistical manifolds and we have given some results for its
induced geometrical objects. Some examples related to these concepts are also presented. Finally, we prove that an invariant lightlike submanifold of indefinite sasakian statistical manifold is an indefinite sasakian statistical manifold.

We hope that, this introductory study will bring a new perspective for researchers and researchers will further work on it focusing on new results not available so far on lightlike geometry

\noindent O\u{g}uzhan Bahad\i r

\noindent Department of Mathematics,

\noindent Faculty of Arts and Sciences,

\noindent Kahramanmaras Sut\c{c}u Imam University

\noindent Kahramanmaras, T\"{u}rkiye

\noindent Email: oguzbaha@gmail.com

\medskip

\end{document}